\documentclass[mathpazo]{cicp}

\usepackage{graphicx}
\usepackage{float}
\usepackage{subfigure}
\usepackage{amssymb}
\usepackage{mathrsfs}
\usepackage{amsmath}
\usepackage{geometry}
\usepackage{titlesec}  
\usepackage{datetime} 
\usepackage[ruled,linesnumbered]{algorithm2e}
\usepackage{color,xcolor}
\usepackage{mlmath}
\usepackage{url}
\usepackage{indentfirst}
\usepackage{color}
\usepackage[colorlinks, linkcolor=red,anchorcolor=blue,citecolor=green]{hyperref}
\usepackage{algpseudocode}

\SetKwRepeat{Do}{do}{while}

\begin{document}

\title{Data-informed Deep Optimization}

\author[Zhang L L et.~al.]{Lulu Zhang\affil{1},  Zhi-Qin John Xu\affil{1,2}\comma\footnotemark,  Yaoyu Zhang\affil{1,2,3}\comma\corrauth}
\address{\affilnum{1}\ School of Mathematical Sciences, Institute of Natural Sciences, Shanghai Jiao Tong University, Shanghai, 200240, China. \\
        \affilnum{2}\ MOE-LSC and Qing Yuan Research Institute, Shanghai Jiao Tong University, Shanghai, 200240, China. \\
          \affilnum{3}\ Shanghai Center for Brain Science and Brain-Inspired Technology, Shanghai, 200031, China}
\emails{{\tt zhangl9661@sjtu.edu.cn} (L.~Zhang), {\tt xuzhiqin@sjtu.edu.cn} (Z.~Xu), {\tt zhyy.sjtu@sjtu.edu.cn} (Y.~Zhang)}

\begin{abstract}

Complex design problems are common in the scientific and industrial fields.
In practice, objective functions or constraints of these problems often do not have explicit formulas, and can be estimated only at a set of sampling points through experiments or simulations. Such optimization problems are especially challenging when  design parameters are high-dimensional due to the curse of dimensionality.
In this work, we propose a data-informed deep optimization (DiDo) approach as follows: first, we use a deep neural network (DNN) classifier to learn the feasible region; second, we sample feasible points based on the DNN classifier for fitting of the objective function; finally, we find optimal points of the DNN-surrogate optimization problem by gradient descent. 
To demonstrate the effectiveness of our DiDo approach, we consider a practical design case in industry, in which our approach yields good solutions using limited size of training data. We further use a 100-dimension toy example to show the effectiveness of our model for higher dimensional problems.
Our results indicate that the DiDo approach empowered by DNN is flexible and promising for solving general high-dimensional design problems in practice.

\end{abstract}


\maketitle

\today


\section{Introduction}

With the development of technology, we are able to study a scientific or industrial problem through massive amount of data \cite{dey2019data}, for example, in neuroscience, modeling a biological neuronal network to meet a set of biological requirements on its dynamical performance \cite{zhang2020dnn}; in industry, optimizing a large set of design parameters to maximize the machine performance while satisfying physical constraints \cite{karen2006hybrid,papalambros2002optimization}. For such practical systems, the dependence between the model or design parameters and the corresponding performance often has no explicit formula \cite{jeong2005efficient,wang2007review}. 
Moreover, constraints in these problems may also be very complex with no explicit formulas. And whether a large set of design parameters is compatible to the constraints can only be examined through experiments or simulations.
Therefore, we often encounter an optimization problem in which both the objective function and constraints are unknown and are available only on a set of data points. For convenience, we call such optimization problems {\it data-informed optimization problem}. Contrary to traditional optimization problems, which are often low-dimensional, can be analytically described and solved by many well-developed algorithms \cite{chong2004introduction,gill2019practical}, it is increasingly important to develop tractable approaches for high-dimensional data-informed optimization problems.

A popular method to solve data-informed optimization problem is to use surrogate models to fit the objective and constraint functions. However, due to the  curse of dimensionality for high dimensional problems, it is often difficult to  obtain a good surrogate model by conventional methods, like polynomial fitting. Empirical and theoretical studies suggest that the DNN approach, trained by gradient-based algorithms, can overcome the curse of dimensionality in fitting high-dimensional functions \cite{ma2019generalization,weinan2020machine}. It has also been observed in practice that the DNNs in general do not overfit even in an overparameterized setting without explicit regularizations \cite{zhang2021understanding}. A series of studies provide possible mechanisms underlying the non-overfitting puzzle of DNNs. For example, frequency principle, both in experiments and theory \cite{xu2020frequency,zhang2019explicitizing, xu2019training, rahaman2019spectral}, shows that DNNs 
prefer to fit training data with low-frequency functions, which often leads to a good generalization performance due to the low frequency dominance in real data. Therefore, DNN serves as an appropriate surrogate model to learn high-dimensional objective and constraint functions from finite samples. We then use a first-order gradient descent to solve the optimization problem. Although this DNN-surrogate optimization problem is non-convex \cite{jain2017non}, we empirically find that the gradient descent can obtain satisfying solutions. We call this DNN-surrogate approach for data-informed optimization problem by {\it data-informed deep optimization} (DiDo). Due to the well-developed program framework, such as Tensorflow and PyTorch \cite{abadi2016tensorflow,paszke2017automatic,paszke2019pytorch} and computation hardwares, such as GPU, the DiDo approach can be easily and efficiently implemented.

The detail of DiDo is briefed as follows. For constraints, since neither the explicit formula nor the number of constraints are completely available, we train a DNN classifier by neural network to learn the feasible region, that is, classifying whether each set of variables satisfy all constraints or not.
To make the DNN classifier more accurate, the training data for the classifier is sampled using an iterative sampling method, that is, adding the data near the boundary of the classifier into training set by Langevin Monte Carlo (LMC)  to re-train the classifier. To improve the efficiency of sampling high-dimensional data within the feasible region for fitting of the objective function by DNN, we sample the training set similarly by LMC based on the well-trained classifier. 
Finally, we find candidates of optimal parameters of the DNN-surrogate optimization problem by traditional algorithm, for example, gradient descent. 

For demonstration, we first consider a specific problem in industry, which is to find the 6-dimensional design parameters for the rotor profile of double screw compressor to maximize the actual flow. 
Without the need of carefully adjusting hyper-parameters, the best actual flow found by our DiDo approach is much better than the result obtained by original hand-craft approach. To illustrate the effectiveness of DiDo approach in higher dimensional problems, we consider a 100-dimension toy example. The optimal value found with our approach is very close to the true optimal value. These results indicate that our DiDo approach can indeed solve high-dimensional data-informed optimization problems.

The rest of paper is organized as follows. Initially, we give a brief preliminary about the notation, DNN and LMC. Followed by the main contents of this paper, our data-informed deep optimization approach, that is, using a deep-based method to solve a type of optimization problem which is different from the traditional ones. Then for demonstrating our DiDo approach, a practical design case in industry and a 100-dimensional toy example are shown in detail. Finally, we make a conclusion and discuss the future work. 

\section{Preliminary}
\subsection{Notation}
In this paper, we use the following notations, see Table \ref{table1}.\\

\begin{table}[htbp]
	\centering  
	\caption{Notation}  
	\label{table1}  
	\begin{tabular}{|c|c|}  
		\hline  
		$x$& scalar component of $\vx$ \\
		\hline  
	
        $\vx$& optimization variable  \\
        \hline  
 
        $d$& dimension of optimization variable  \\
        \hline  
 
        $[n]$& index set $\{1,2,...,n\}$ \\
        \hline  
 
        $f(\vx)$& objective function  \\
        \hline  
 
        $\Omega$& feasible region  determined by the considered problem \\
        \hline  
 
        $\partial \Omega$& the true boundary of implicit feasible region \\
        \hline  
 
        ${I}_{\Omega}(\vx)$& indicator function of the region $\Omega$, i.e., if $\vx \in \Omega$,${I}_{\Omega}(\vx)$=1; otherwise,${I}_{\Omega}(\vx)$=0 \\
        \hline  
 
        $D_{\mathrm {obj}}=\{(\vx_i;f(\vx_i)\}_{i=1}^{n_o}$& training set for DNN fitting \\
        \hline  
 
        $D_{\mathrm {c}}=\{(\vx_i;{I}_{\Omega}(\vx_i))\}_{i=1}^{n_c}$& training set for DNN classifier \\
        \hline  
 
        $f_{\vtheta_{\mathrm{o}}}(\vx)$& DNN surrogate model for objective function\\
        \hline  
 
        $f_{\vtheta_{\mathrm{c}}}(\vx)$& DNN classifier neural network for feasible region\\
        \hline  
 
		\hline
	\end{tabular}
\end{table}


\subsection{DNN}
The general setup for a DNN is reviewed as follows. A fully connected DNN of $H$ layers is denoted by 
\begin{align*}
&f_{\vtheta}(\vx) = \vW^{[H-1]} \sigma \circ (\vW^{[H-2]} \sigma \circ(\cdots(\vW^{[1]} \sigma \circ(\vW^{[0]} \sigma \circ + \vb^{[0]})+\vb^{[1]})\cdots)+\vb^{[H-2]})+\vb^{[H-1]},
\end{align*}
where  $\vx \in \mathbb{R}^d$, $\vW^{[l]} \in \mathbb{R}^{m_{l-1}\times m_l}$, $\vb^{[l]} \in \mathbb{R}^{m_{l-1}}$, $m_0=d$, $m_H=1$, $\sigma$ is the activation function and $``\circ"$ means entry-wise operation.
The set of parameters for DNN is denoted by
\begin{equation*}
\vtheta=(\vW^{[0]}, \vW^{[1]},\cdots, \vW^{[H-1]},\vb^{[0]},\vb^{[1]},\cdots,\vb^{[H-1]}),
\end{equation*}

For the regression problem of fitting a training set $\{(x_i,y_i)\}_{i=1}^{n}$, where $x_{i}\in\mathbb{R}^{d}$ and $y_{i} \in \mathbb{R}$ for each $i$, the commonly used loss functions are  mean-square error (MSE), that is,
\begin{equation*}
L(\vtheta)=\frac{1}{n}\sum_{i=1}^{n}(f_{\vtheta}(\vx_i)-y_i)^2,
\end{equation*}
and root-mean-square error (RMSE), that is,
\begin{equation*}
L(\vtheta)=\sqrt{\frac{1}{n}\sum_{i=1}^{n}(f_{\vtheta}(\vx_i)-y_i)^2}.
\end{equation*} 
We use MSE when training DNN to fit the objective function. 
and use RMSE to measure the training error and test error of the DNN.

For the  classification problem of fitting a training set $\{(x_i,q_i)\}_{i=1}^{n}$, where $x_{i}\in\mathbb{R}^{d}$ and $q_{i} \in \{0,1\}$ for each $i$, the loss function we used is binary-cross-entropy (BCE), that is,
\begin{equation*}
L(\vtheta)=-\frac{1}{n}\sum_{i=1}^n[q_ilogf_{\vtheta}(\vx_i)+(1-q_i)log(1-f_{\vtheta}(\vx_i))]
\end{equation*}

In cases given in this paper, the activation function is fixed to GELU function.
GELU is a smooth non-saturating activation function that can alleviate gradient vanishing. Empirically, a GELU DNN is efficient to train and generalizes well for smooth problems we considered.
 Note that one can also consider other smooth non-saturating activation functions like Swish, ELU and SELU to achieve similar training and generalization performance.

 During the training of neural network, the parameters of the DNN in each epoch  are updated by a gradient-based optimization algorithm, e.g. gradient descent (GD), stochastic gradient descent (SGD) or Adam.
To speed up the training process,  we update the parameters of DNN using Adam \cite{kingma2014adam}.

\subsection{Langevin Monte Carlo (LMC)}
There are many mature methods to sample data from a desired probability distribution, such as Markov chain Monte Carlo, Metropolis-Hastings, Hamiltonian Monto Carlo and Split Monte Carlo. 
For convenience, in our experiments, we use overdamped Langevin Monte Carlo (LMC).

LMC is a common method to sample data  following a Boltzmann distribution \cite{roberts1996exponential,dalalyan2017theoretical, durmus2016sampling, dalalyan2019user}. This method is based on evolving a stochastic differential equation (SDE), that is, 
\begin{equation*}
d\vx=-\nabla E(\vx) dt + \sqrt{\frac{2}{\beta}}dW,
\end{equation*}
where  $\beta$ is positive hyperparamter and $W$ is the Brownian motion.
The steady-state distribution of this SDE is proportional to $e^{-\beta E(\vx)}$ and it satisfies the detailed balance condition. 
The set of long-time solution of the SDE follows the Boltzmann distribution $\sim e^{-\beta E(\vx)}$. 
We use the first-order  Euler-Maruyama scheme to solve the SDE, i.e. each data point $\vx$ is updated according to $\vx^{t+1}=\vx^{t}-\alpha \nabla E(\vx^{t}) + \sqrt{\frac{2 \alpha}{\beta}}\xi^{t}$, where $\xi^t \sim N(0_d,I_d)$.

Note that we choose appropriate energy function $E(\vx)$ for different tasks. In our experiment, we use  $E(\vx)=(f_{\vtheta_{\mathrm{c}}}(\vx)-0.5)^2$ to sample data  concentrated around the boundary of the predicted-feasible region for DNN classifier $f_{\vtheta_{\mathrm{c}}}(\vx)$ and use $E(\vx)=(f_{\vtheta_{\mathrm{c}}}(\vx)-1)^2$ to efficiently sample more feasible data, where $f_{\vtheta_{\mathrm{c}}}(\vx)$ denotes the DNN classifier.

A simple version of LMC method is shown in algorithm \ref{alg:LMC}.

\begin{algorithm}
\caption{Langevin Monte Carlo (LMC)\label{alg:LMC}}
\KwData{$T$: total iteration steps; energy function $E(\vx)$; step length $\alpha$; positive $\beta$; $X^0=\{\vx^0_i\}_{i=1}^n$:initial data set.}
\KwResult{$X^T$}
\For{$t = 0$ \KwTo $T$}{$\vx^{t+1}_i=\vx^t_i-\alpha  \nabla E(\vx^t_i) + \sqrt{\frac{2 \alpha}{\beta}}\xi^t_i$, where $\xi^t_i \sim N(0_d,I_d), i\in [n]$ \;}
Get $X^{\rm T}=\{\vx^{\rm T}_i\}_{i=1}^n$ and $X^{\rm T} \sim e^{-\beta E(\vx)}$
\end{algorithm}

 

\section{Data-informed deep optimization}
In this section, we introduce the framework of Data-informed deep optimization (DiDo) approach for solving high-dimensional optimization problems, in which the objective function and the constraints are only available through samples without explicit formula. The DiDo approach shows an indispensable value beyond the tradition optimization approach in high-dimensional data-informed problems.

\subsection{Data-informed problem formulation}

The data-informed optimization problem is formulated as follows. 

\noindent \textbf{Data-informed optimization problem:}
\begin{equation}
\mathop{\min}\limits_{\vx \in \Omega}\quad f(\vx),
\label{data-informed problem}
\end{equation}
where objective function $f(\vx)$ and feasible region $\Omega$ are implicit which can only be evaluated through simulation at certain sampling points. In practice,  $\Omega$ is often defined by a series of implicit constrains as  $\Omega=\{\vx|f_i(\vx)\leq 0,i=1,2,...,L\}$.

\subsection{Deep optimization approach}
\begin{figure}[!htb]
    \centering
    \includegraphics[scale=0.05]{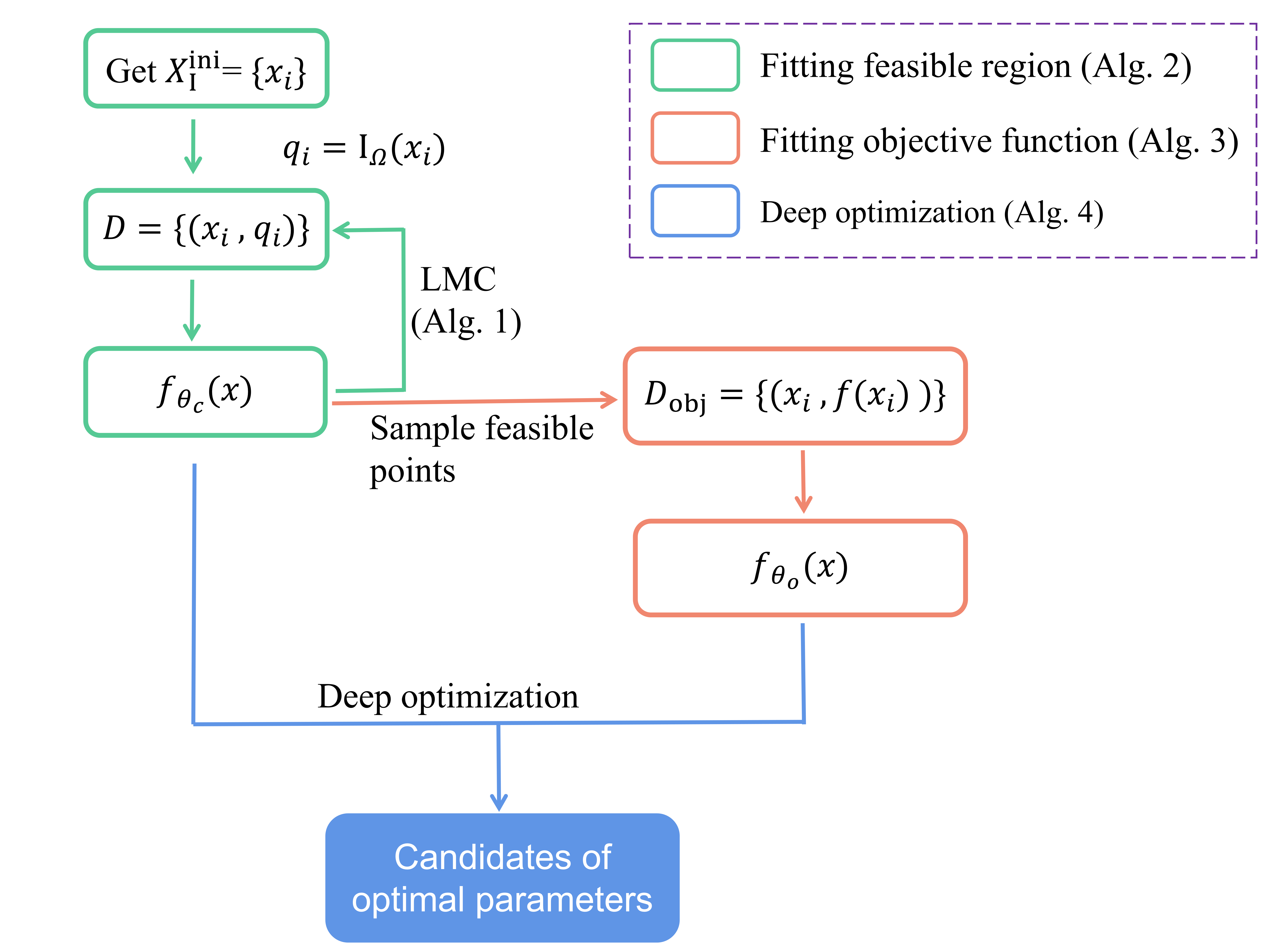}
  
    \caption{ The flow chart of the DiDo approach}
    \label{fig: flowchart}
\end{figure}

We propose a deep optimization approach to solve the data-informed optimization problem (see Fig.  \ref{fig: flowchart} for a flow chart).
Our general idea to deal with a problem with implicit objective function and feasible region is to fit them by DNN surrogates from data. Then we can optimize this problem with common gradient-based method. Note that the training of the objective function relies on an explicit and accurate surrogate of feasible region for generating high quality training samples well covering the whole feasible region. Therefore, we first train a  DNN classifier $f_{\vtheta_{\mathrm{c}}}(\vx)$  through an iterative process from an initial sample set.  Then we use $f_{\vtheta_{\mathrm{c}}}(\vx)$  to generate random samples from feasible region and evaluate the corresponding values of the objective function. Then we build the training data set $D_{\mathrm {obj}}=\{\vx_i,f(\vx_i)\}_{i=0}^{n_o}$ by simulation, whose inputs are sampled randomly from feasible region based on  $f_{\vtheta_{\mathrm{c}}}(\vx)$. Through fitting $D_{\mathrm {obj}}$, we obtain a DNN  $f_{\vtheta_{\mathrm{o}}}(\vx)$ as a surrogate of the objective function $f(\vx)$.
Finally, by optimizing $f_{\vtheta_{\mathrm{o}}}(\vx)$ with surrogate constrains $f_{\vtheta_{\mathrm{c}}}(\vx) \geq 0.5$ ($f_{\vtheta_{\mathrm{c}}}(\vx) = 0.5$ is regarded as surrogate boundary) , we can get candidates of the optimal parameters of the problem (\ref{data-informed problem}), which should be close to the true optimal parameters of the problem.

\subsubsection{Fitting feasible region}

Generally, without an explicit feasible region, it is difficult to generate well distributed feasible training samples especially in a high-dimensional problem. With a blind sampling, training samples are likely
far from the decision boundary, i.e., boundary of the feasible region, resulting in an inaccurate fitting of the DNN classifier. 
To overcome this difficulty in our  deep optimization  approach, we propose an iterative method which adds new samples around the boundary of current DNN classifier and retrain it at each iteration. Using this approach, we can efficiently obtain a good DNN classifier $f_{\vtheta_{\mathrm{c}}}(\vx) \in [0,1]$ through several rounds of iteration.

Initially, we uniformly sample $X^{\rm ini}_{\rm I}$ in a selected region $B$ based on the prior knowledge of the considered problem and train the classifier $f_{\vtheta_{\mathrm{c}}^{(0)}}(\vx)$ by $D^{(0)}=\{(\vx_i,{I}_{\Omega}(\vx_i))|{\vx}_i \in X^{\rm ini}_{\rm I}\}$. Empirically balancing the feasible and infeasible points benefits the performance of the classifier. 
Note that  many problems whose optimal parameters close to the boundary of the feasible region require highly accurate DNN classifier (see example in Fig.  \ref{fig:6para-optimize-boundary-divergence}).
We propose an iterative method to efficiently improve the accuracy of classifier $f_{\vtheta_{\mathrm{c}}^{(t)}}(\vx)$  at each iteration step $t$. 
For classification problem, generally, the points close to the decision boundary is of crucial importance to determine the classifier, e.g., support vectors for support vector machine (SVM). 
Therefore,  at each iteration step, we add new training data sampled near the decision boundary of classifier $f_{\vtheta_{\mathrm{c}}^{(t)}}(\vx)$ by LMC method (see algorithm \ref{alg:LMC} for details) and train a new classifier $f_{\vtheta_{\mathrm{c}}^{(t+1)}}(\vx)$ initialized by ${\vtheta_{\mathrm{c}}^{(t)}}$.  



\begin{figure}[!htb]
    \centering
    {
    \includegraphics[scale=0.06]{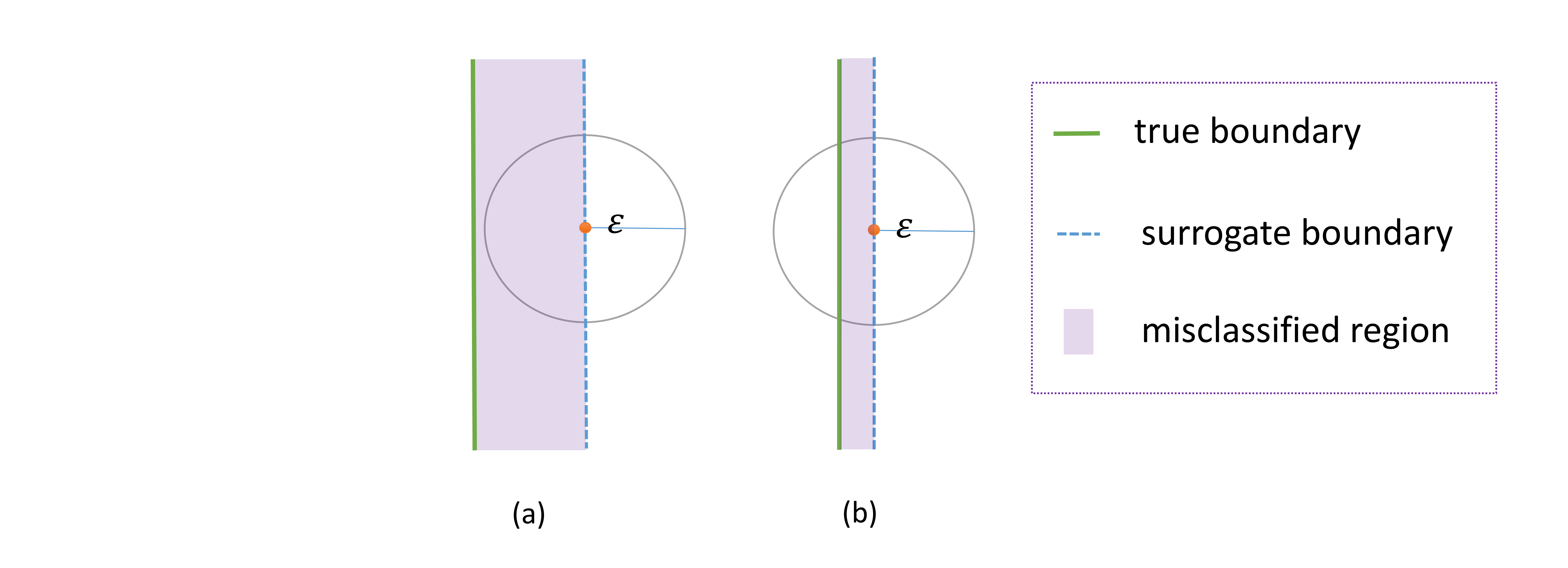}
    }
    \caption{schematic diagram. (a) For red point on the surrogate boundary, the distance  to the true boundary is larger than $\epsilon$, and the  prediction accuracy is roughly $50\%$; (b) for red point on the surrogate boundary, the distance  to the true boundary is much smaller than $\epsilon$, and the  prediction accuracy is close to $1$. }
    \label{fig:schematicdiagram}
\end{figure}

For a  stopping criterion, it is crucial to determine whether the surrogate boundary is close to the true boundary, e.g.,  their ``mean distance'' is smaller than  certain tolerance $\epsilon$.
Intuitively, for any point on the surrogate boundary, if its distance to the true boundary is larger than the  $\epsilon$, then the prediction accuracy of the DNN classifier in the $\epsilon$-neighborhood of this point is roughly $50\%$ (see Fig. \ref{fig:schematicdiagram}(a) for illustration); otherwise, if the distance is much smaller than $\epsilon$,  then the  prediction accuracy in the $\epsilon$-neighborhood should be close to $1$ (see Fig. \ref{fig:schematicdiagram}(b) for illustration).
Therefore, 
we sample some points close to the surrogate boundary by LMC method (see algorithm \ref{alg:LMC} for details) and perturbed them by Gaussian noise of covariance matrix $\sigma^2I_d$, where $\sigma$ is roughly $\epsilon$ due to concentration in the equator \cite{blum2020foundations}. 
 When the predicted accuracy of the classifier on these points is higher than a expected  value, say $95\%$, we stop the iteration.

The detail of our iterative training algorithm is shown in algorithm \ref{alg:iterative method}.

\begin{algorithm}
\caption{Iterative training algorithm\label{alg:iterative method}} 
\KwData{$B$: region of initial sampling depending on the considered problem; $n_0$: initial sample size; $n_{1}$: adding sample size at each iteration; $E_t(\vx) = {(f_{\vtheta_{\mathrm{c}}^{(t)}}(\vx)-0.5)}^2$: energy function used in LMC method; $\sigma$:  standard deviation of noise term; $\beta$: positive hyperparameter used in LMC method.}
\KwResult{Good classifier: $f_{\vtheta_{\mathrm{c}}}(\vx)$}
Uniformly sample $X^{\rm ini}_{\rm I}$ in $B$\;
Define $D^{(0)}=\{(\vx_i,{I}_{\Omega}(\vx_i))|{\vx}_i \in X^{\rm ini}_{\rm I},i\in[n_0]\}$\;
Define  $t=0$\;

\Do {}{
Train $\vtheta_{\mathrm{c}}^{(t)}$ of $f_{\vtheta_{\mathrm{c}}^{(t)}}(\vx)$ by $D^{(t)}$ with Adam\;
Use LMC method with proper initialization  and 
  $E_t(\vx)$ to sample $n_1$ data following distribution $\sim e^{-\beta E_t(\vx)}$ and obtain input set $X^{(t)}_{\rm I}$\;
      Perturbation: $X_{P}^{(t)}=\{\vx + \vxi|\vx \in X^{(t)}_{\rm I}, {\vxi} \sim {\mathcal N}(0_d,{\sigma}^2I_d)\}$\;
        Evaluate the classification accuracy $acc$  of $f_{\vtheta_{\mathrm{c}}^{(t)}}(\vx)$ on $X_{P}^{(t)}$ \;
       \If{$acc \geq  95 \%$}{
       $f_{\vtheta_{\mathrm{c}}}(\vx) \leftarrow f_{\vtheta_{\mathrm{c}}^{(t)}}(\vx)$\;
            break\;}
     Evaluate  $X^{(t)}_{\rm I}$ and add to the training data $ D^{(t+1)} =\{(\vx_i,{I}_{\Omega}(\vx_i))|\vx_i \in X^{(t)}_{\rm I} \}\cup D^{(t)}$\;
     
      Update $t \leftarrow t+1$ \;
    }

\end{algorithm}

\subsubsection{Fitting objective function}

For a  high-dimensional large-scale problem, with implicit boundary, it is difficult to efficiently sample  diverse training data.  
However, with explicit classifier obtained above, we can use LMC with energy function  $E(\vx) = {(f_{\vtheta_{\mathrm{c}}}(\vx)-1)}^2$ to generate high quality training samples $D_{\mathrm {obj}}=\{(\vx_i,f(\vx_i)\}$ well covering the feasible region of considered problem.
By training the DNN by  $D_{\mathrm {obj}}$, we can get the DNN surrogate $f_{\vtheta_{\mathrm{o}}}(\vx)$ of the objective function.

The detail of  fitting objective function is shown in  algorithm \ref{alg:sample method}.




\begin{algorithm}
\caption{Fitting of objective function} \label{alg:sample method} 
\KwData{classifier $f_{\vtheta_{\mathrm{c}}}(\vx)$; a non-empty feasible set of   $f_{\vtheta_{\mathrm{c}}}(\vx)$: $X^{\rm ini}_{\rm S}$; a large enough number $n_t$; energy function $E(\vx) = {(f_{\vtheta_{\mathrm{c}}}(\vx)-1)}^2$.}
\KwResult{DNN surrogate $f_{\vtheta_{\mathrm {o}}}(\vx)$}
Generate initial points for LMC method: $X_0 = \{\vx_i| \vx_i \in X^{\rm ini}_{\rm  S}, i\in [n_t]  \}$\; 
Use LMC method with $E(\vx)$ to sample data following distribution $\sim e^{-\beta E(\vx)}$ and obtain $X^{\rm T}_{\rm S}$\;
Select data in the feasible region of real system $\Omega$: $X_{\mathrm{o}}=\Omega \cap X^{\rm T}_{\rm S}$\;
Obtain training data for the objective function: $D_{\mathrm {obj}}=\{(\vx_i,f(\vx_i))|\vx_i \in X_{\mathrm{o}},i\in [n_{t'}]\}$\;
Train $\vtheta_{\mathrm {o}}$ of $f_{\vtheta_{\mathrm {o}}}(\vx)$ by $D_{\mathrm {obj}}$ with Adam.
\end{algorithm}




  

\subsection{Deep optimization}

Based on the  accurate DNN surrogate models of constraints and objective function obtained above, the data-informed  optimization problem (\ref{data-informed problem}) turns to be the following explicit optimization problem,
\begin{align}
&\min_{\vx}\quad f_{\vtheta_{\mathrm{o}}}(\vx) \notag\\
&\hspace{0.1cm}\mathrm{s.t.}\hspace{0.3cm}0.5-f_{\vtheta_{\mathrm{c}}}(\vx) \leq 0, 
\label{subprob1}
\end{align}
where $0.5$ is the threshold of the DNN classifier $f_{\vtheta_{\mathrm{c}}}(\vx) \in [0,1]$ for prediction. 

The problem (\ref{subprob1}) is a conventional optimization problem with constraints. To solve it, we first rewrite it as an unconstrained problem, making the inequality constraint implicit in the objective
\begin{equation*}
\min_{\vx}\quad f_{\vtheta_{\mathrm{o}}}(\vx) + I_{-}(0.5-f_{\vtheta_{\mathrm{c}}}(\vx)), \\
\end{equation*}
where $I_{-}:\mathbb{R} \mapsto \mathbb{R}$ is the indicator function for the non-positive real number,
\begin{equation*}
 I_{-}(u)= \begin{cases}
0,&\text{if} \quad u \leq 0;\\
\infty,&\text{if} \quad u > 0.\
\end{cases}
\end{equation*}
However, the indicator function $I_{-}$ is not differentiable. We approximate the indicator function $I_{-}$ by a ``soft'' function. 
 For example, we use the interior-point method. The basic idea of interior-point method is to  approximate the indicator function $I_{-}(u)$ by the barrier function and a common barrier function is logarithmic barrier, $-(\frac{1}{t})log(-u)$, where $t>0$ is a hyperparameter that sets the accuracy of the approximation \cite{boyd2004convex}. 
 
 Substituting $I_{-}(u)$ with $-\frac{1}{t}log(-u)$ gives the approximation
 \begin{equation}
\min_{\vx}\quad f_{\vtheta_{\mathrm{o}}}(\vx) -(\frac{1}{t})log(-(0.5-f_{\vtheta_{\mathrm{c}}}(\vx))). \\
\label{subprob2}
\end{equation}
To solve  problem (\ref{subprob2}), 
we use gradient descent (GD) for convenience. Although simple, we find that GD is often an effective optimization algorithm in DiDo.

The  deep optimization is concluded in algorithm \ref{alg:Deepoptimization}.

\begin{algorithm}
\caption{Deep optimization \label{alg:Deepoptimization}} 
\KwData{$f_{\vtheta_{\mathrm{c}}}(\vx)$: well-trained DNN classifier; $f_{\vtheta_{\mathrm{o}}}(\vx)$: DNN surrogate model for fitting objective function.}
\KwResult{Candidates of  optimal parameters}
Substitute  $f_{\vtheta_{\mathrm{o}}}(\vx)$ and $f_{\vtheta_{\mathrm{c}}}(\vx)$ into problem (\ref{subprob2})\;
Solve problem  (\ref{subprob2})  by gradient-descent-based optimization algorithms, such as gradient descent (GD)\;
Get candidates of  optimal parameters of the problem (\ref{data-informed problem}).
\end{algorithm}

  
   

The proposed methodology gives a schematic process to search for candidates of optimal parameters (see Fig. \ref{fig: flowchart}) for high dimensional optimization problem with implicit feasible region and objective function.  
As we will show, it is well suitable for data-driven inferences using deep neural networks which can efficiently differentiate.

Remark that even when we can analytically characterize the feasible region by a set of equations, we  can also train a DNN surrogate to represent the feasible region. In such case, our approach can still bring benefits, for example, using DNN classifier can soft the boundary of the feasible region and we can easily determine the normal vector of the boundary.

\section{Optimal rotor profile}
In this section, we apply the DiDo approach to solve an engineering design problem to show its effectiveness.  
\subsection{Problem description}
 Screw compressor is widely used in refrigeration, mining, petrochemical and other industries because of its high reliability, good power balance, less leakage and high efficiency. As the core component of twin-screw compressor, optimizing the design of rotor profile would vastly benefit the mechanical performance of the screw compressor. The rotor profile is smoothly connected by several arcs and arc envelopes together. Empirically, we can parameterize the rotor profile by 6 parameters, $\vx=[r, r_{\mathrm{3}}, r_{\mathrm{o}}, r_{\mathrm{o2}}, u_{\mathrm{1}}, R] \in \mathbb{R}^6$, where $r$, $r_{\mathrm{3}}$, $r_{\mathrm{o}}$, $r_{\mathrm{o2}}$, $R$ are  radius of the arc and $u_{\mathrm{1}}$ is an angle \cite{xing2000screw,wu2004theoretical}. Then, the optimization of the rotor profile becomes an optimization problem w.r.t. the 6 parameters. 
 
 In our example, the performance of a design parameter set, consisting of the 6 design parameters, is measured by the actual flow of the rotor, which is an important performance indicator for large compressor, through computational fluid dynamics simulation. Our goal is to find a rotor profile that can maximize the actual flow.

Remark that not all parameters in $\mathbb{R}^6$ are feasible for the design.   They should satisfy a set of implicit constraints related to geometrical properties of the rotor.
 Therefore, both the objective and the constraint functions are data-informed, i.e., they are only available on a set of data points through simulation. 
In the following, we demonstrate the effectiveness of our DiDo approach on this problem. 
 
\subsection{Feasible region learned by a DNN classifier}
In this example, we first use the iterative training algorithm in  algorithm \ref{alg:iterative method} to train the DNN classifier $f_{\vtheta_{\mathrm{c}}}(\vx)$, which is a fully connected DNN with hidden layer sizes 800-600-400-200 equipped with a sigmoid function at the output layer. Without loss of generality, we choose $0.5$ as threshold to determine the surrogate feasible region, i.e., $f_{\vtheta_{\mathrm{c}}}(\vx) \geq 0.5$.

  Remark that, we carefully choose the initial sample region $B$, such that  the number of feasible points and  non-feasible points are balanced in the initial training data. For the effectiveness of DNN training,  we normalize each parameter to a mean zero and  variance one input variable.

We set initial sample size $n_0=8000$ and  we set $n_1=5000$ samples in each iteration.  With algorithm \ref{alg:iterative method}, we can obtain a well-trained classifier  $f_{\vtheta_{\mathrm{c}}}(\vx)$.

\begin{figure}[!htb]
    \centering
    \subfigure[]{
    \includegraphics[scale=0.15]{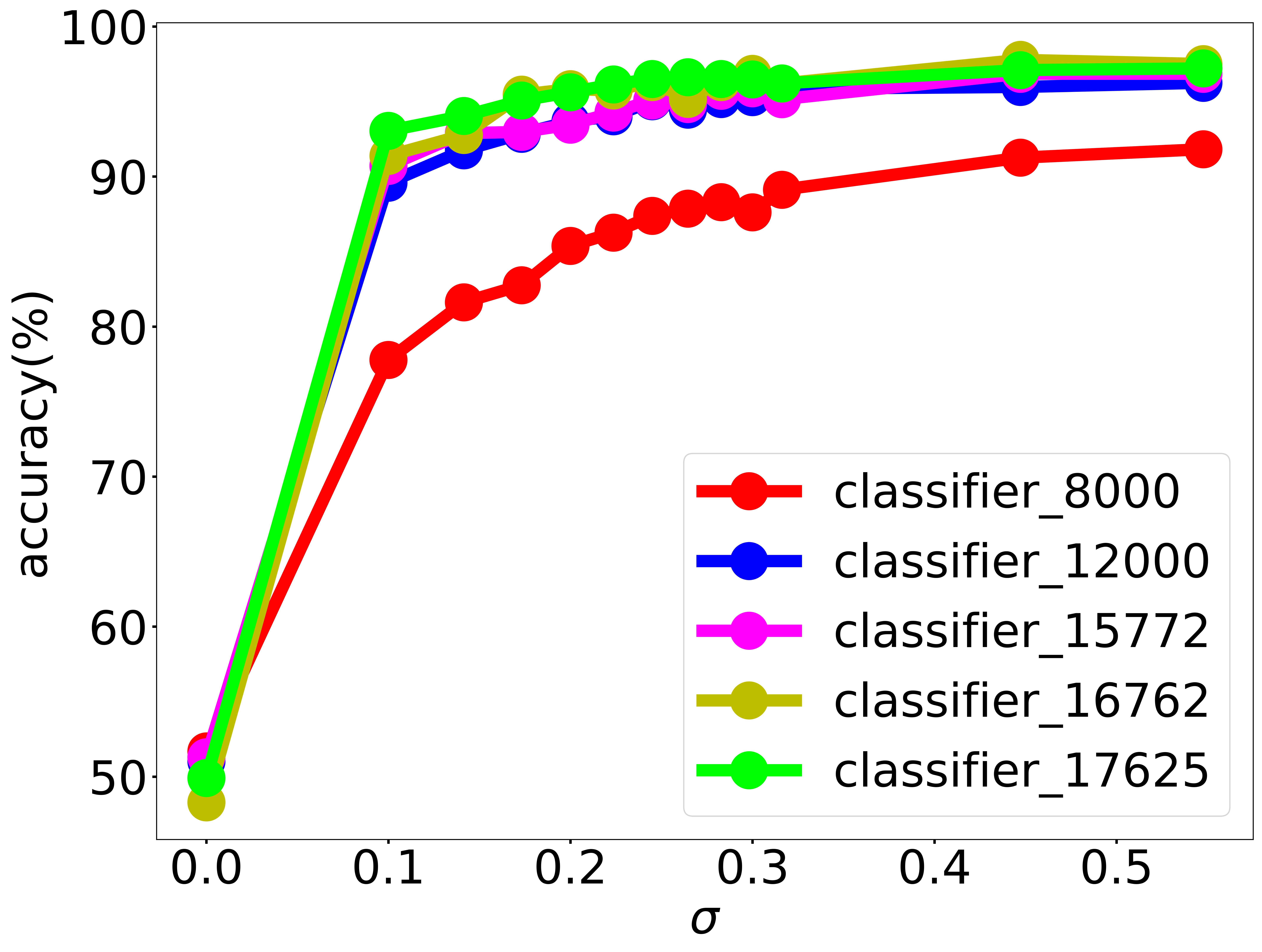}
    }
    \quad
    \subfigure[]{
    \includegraphics[scale=0.15]{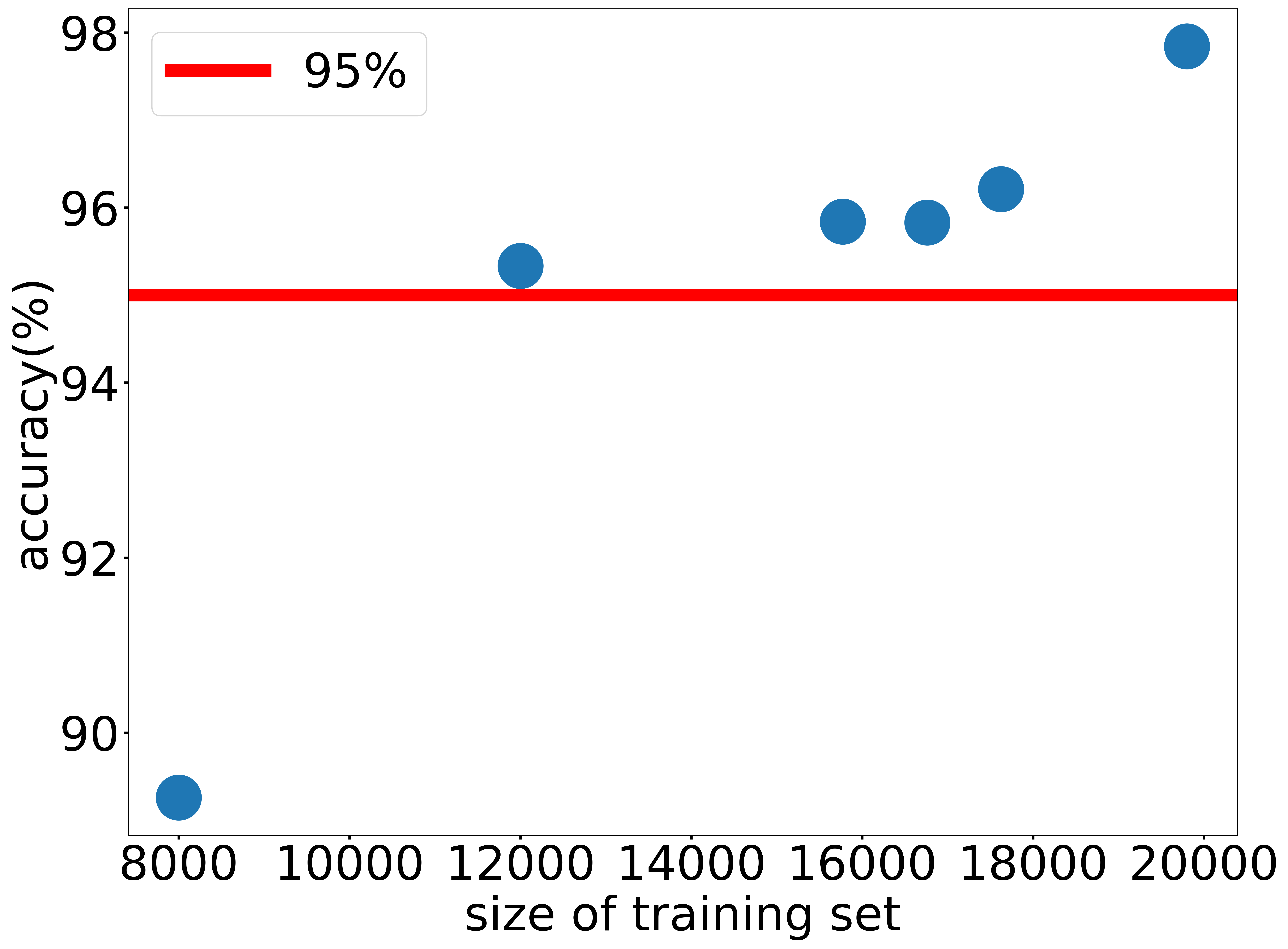}
    }
    
    \caption{ (a) Classification accuracy of the DNN classifier on the perturbed terms during iteration. Note that, at each iteration $t$, we apply an extra constraint  $|f_{\vtheta_{\mathrm{c}}}^{(t)}(\vx_i)-0.5|\leq 0.1\}$ to the points sampled by LMC.
In the two figures, label accuracy means classification accuracy after perturbation.  As we add more data , the magnitude of the perturbed term when classifier accuracy on perturbed term achieve  $100\%$ gets smaller, which means the performance of classifier is better. (b) Classification accuracy of the DNN classifier on the fixed standard deviation of the perturbed terms, where  variance  $\sigma^2 = 0.1$. The classification accuracy is getting better as we update the DNN classifier.}
    \label{fig: 6para-classifier-performance}
\end{figure}

To show effectiveness of the iterative training, we show the accuracy of the DNN classifier on the  samples at surrogate boundary at each iteration with Gaussian noise perturbation during the iteration. As shown in Fig. \ref{fig: 6para-classifier-performance}(a), for each curve, which is the accuracy of the classifier w.r.t. different noise standard deviation, as the perturbation noise increases, the accuracy increases. This indicates that the classifier is more accurate on the samples that deviate more from the boundary, which provides a rationale for our iterative training algorithm focusing on training the boundary. Compared with different iterations, indicated by different colored curves, as the iteration proceeds accompanied by the increasing of training samples, the classifier is improved. For example, as shown in Fig. \ref{fig: 6para-classifier-performance} (b), considering a fixed noise 
with variance $0.1$, the accuracy of the classifier almost monotonically increases as the size of the training set.

\subsection{Objective function learned by a DNN model}

We use a DNN surrogate  to fit the objective function, i.e. a mapping from a designed rotor profile to the actual flow.
 By algorithm \ref{alg:sample method}, we use  the classifier $f_{\vtheta_{\mathrm{c}}}(\vx)$ obtained above to generate a training set $D_{\mathrm {obj}}$ of size $n_t = 500$  
and train a GELU-DNN $f_{\vtheta_{\mathrm{o}}}(\vx)$ of hidden layer size 1024-512-256-128.
The test accuracy of the DNN $f_{\vtheta_{\mathrm{o}}}(\vx)$  is evaluated on a test data set consisting of 2000 samples. 

\begin{figure}[!htb]
    \centering
    \includegraphics[scale=0.15]{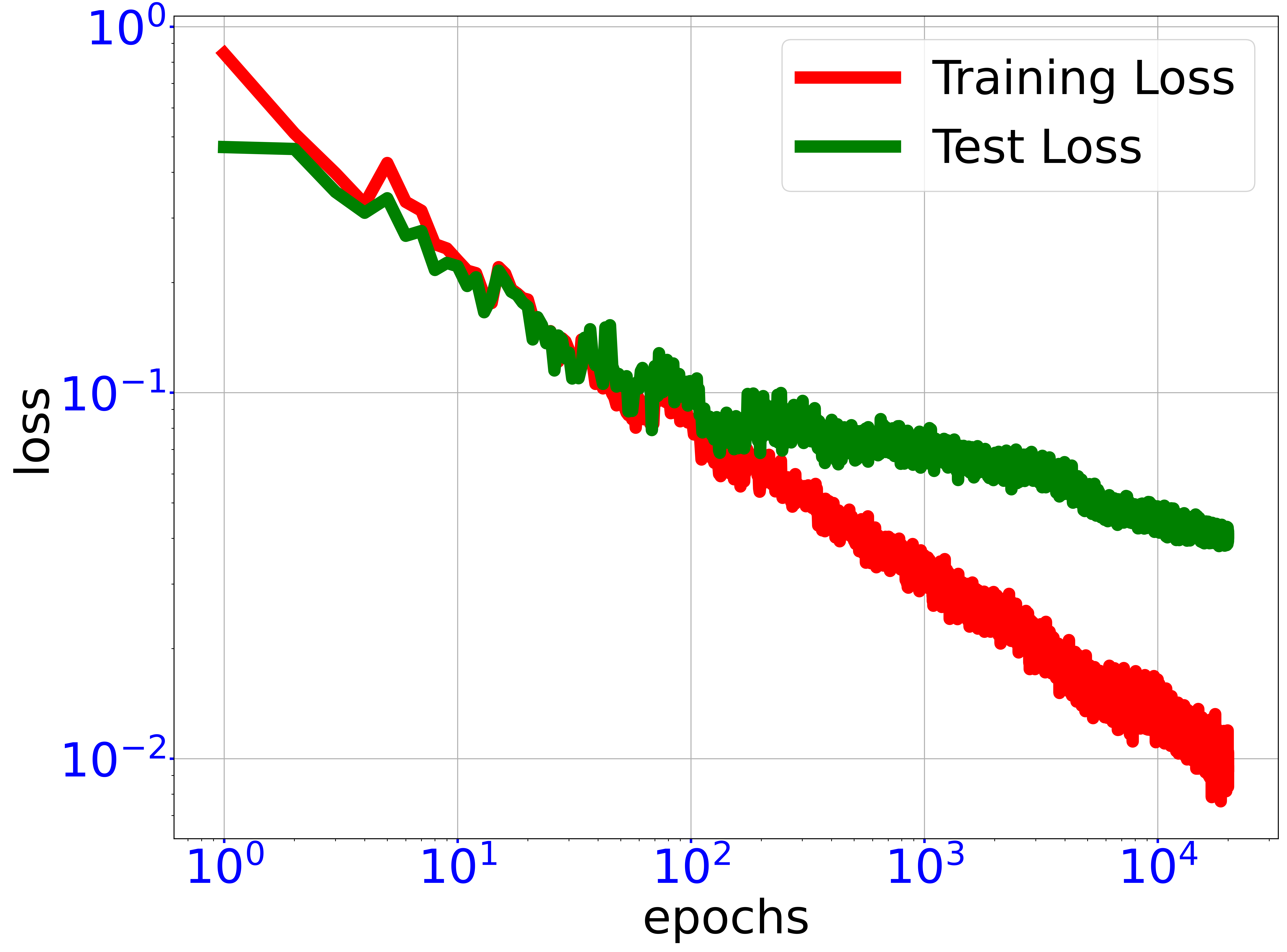}
    \caption{Training loss and test loss during training  DNN surrogate $f_{\vtheta_{\mathrm{o}}}(\vx)$ of optimal rotor profile problem.}
    \label{fig:rotorDNNloss2}
\end{figure}

As shown in  Fig. \ref{fig:rotorDNNloss2}, after training, 
 the normalized RMSE training error is $\sim 0.01$ whereas the normalized RMSE test error is $\sim 0.04$.

\subsection{Deep optimization}

Then we solve the problem with data-informed deep optimization approach in algorithm \ref{alg:Deepoptimization} using $f_{\vtheta_{\mathrm{o}}}(\vx)$ and $f_{\vtheta_{\mathrm{c}}}(\vx)$.
 
\begin{figure}[H]
    \centering
    \includegraphics[scale=0.15]{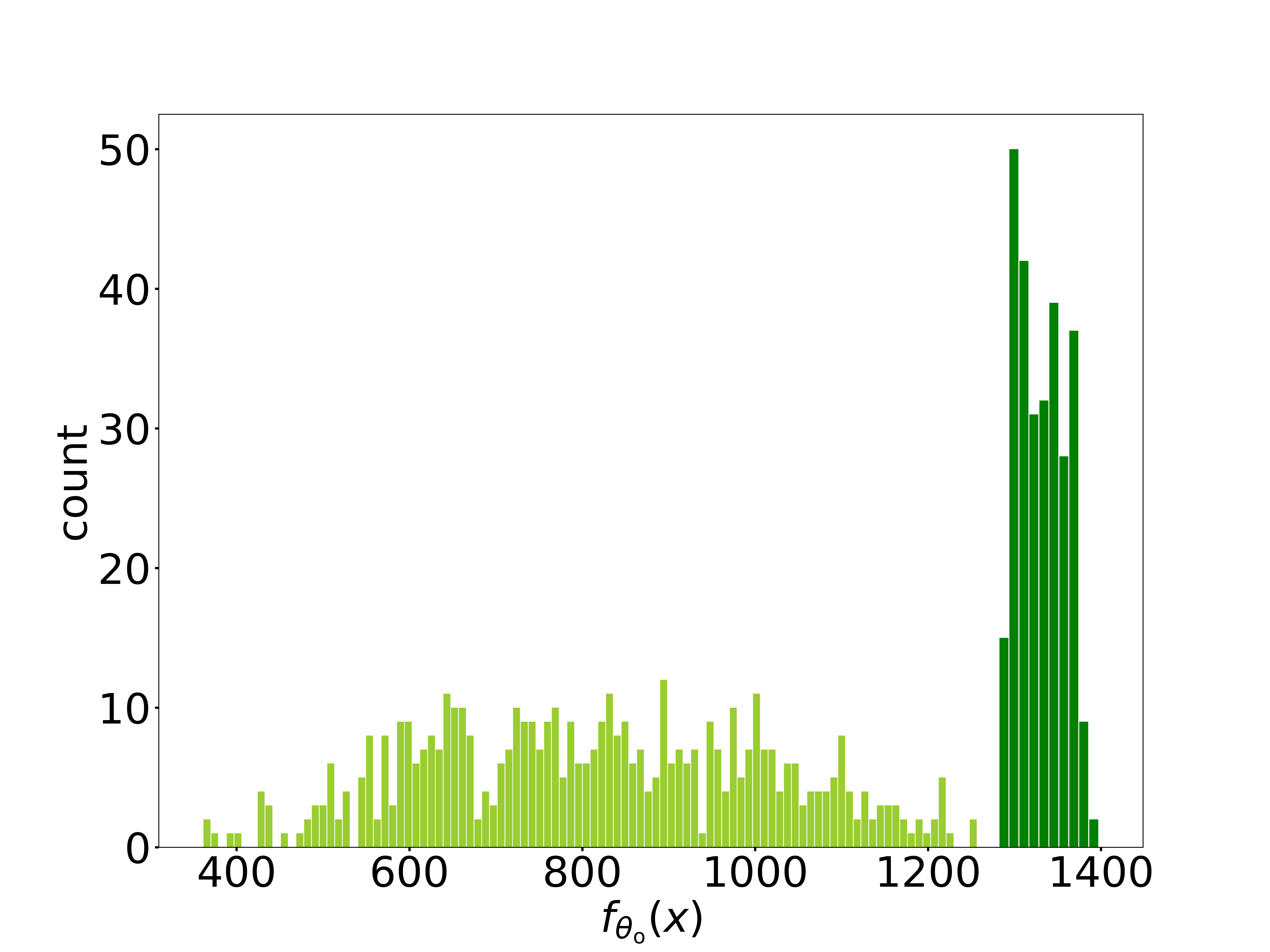}
    \caption{ The distribution of the simulation actual flow value  on  sampled data used for training DNN surrogate and the candidates of optimal parameters obtained finally. The light green bars correspond to the training samples and  the dark green bars correspond to the candidates of optimal parameters. }
    \label{fig:6para-optimize-result}
\end{figure}

The optimal of this optimization problem may be not unique and there could be multiple local minima. 
Therefore, we solve the problem by gradient descent with various initial points to search for a global minimum. For visualization, in Fig. \ref{fig:6para-optimize-result}, we show the distribution of the actual flow
of the training samples used for learning DNN surrogate and a set of true feasible candidates of optimal profile parameters. Note that the maximal actual flow of training samples  approximates $1256$. 
After solving the optimization problem, we obtain a set of candidates of optimal profile parameters. 
Then we examine whether those parameters are in true feasible region with simulator and calculate the actual flow on these feasible  designed rotor parameters with CFD simulator. The best actual flow we achieved is roughly $1400$, which is better than those obtained by manually tuning parameters and the maximal actual flow of training samples $1256$. The candidates of optimal profile parameters outperform the training samples in the sense of the actual flow. Most of the actual flow of the candidates of optimal profile parameters are larger than $1340$. 


%

\begin{figure}[H]
    \centering
    {
    \includegraphics[scale=0.15]{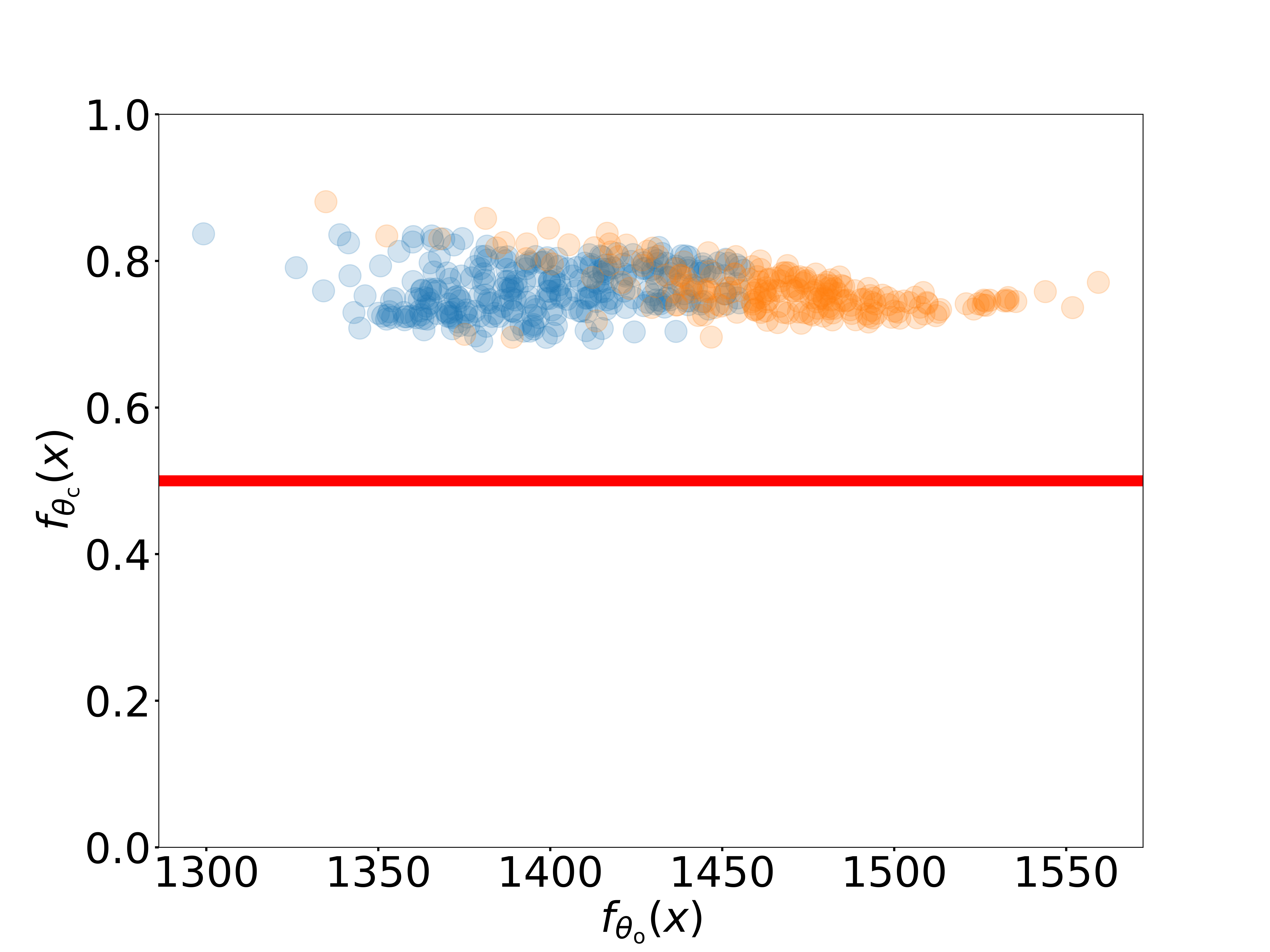}
    }
    \caption{ The classifier value $f_{\theta_{c}(\vx)}$ and the actual flow  predicted by DNN surrogate $f_{\theta_{o}(\vx)}$ on these candidates of optimal parameters. The red solid line is corresponding to the probability 0.5. Both blue and yellow dots are feasible predicted by DNN,
    both above the solid red line.
    However the yellow points are outside the true boundary.
    }
    \label{fig:6para-optimize-boundary-divergence}
\end{figure}

Further more, it is interesting to analyze the candidates of optimal 
parameters obtained using our DiDo approach. We analyze the distance between the candidates of optimal parameters and the boundary of the feasible region by computing the probability predicted by the classifier $f_{\vtheta_{\mathrm{c}}}(\vx)$. As shown in Fig. \ref{fig:6para-optimize-boundary-divergence}, 
each  point corresponds to a candidates of optimal parameter
and the  $f_{\vtheta_{\mathrm{c}}}(\vx)$ of obtained candidates of optimal parameters    
significantly deviate from 1,  i.e., most candidates of optimal parameters with different actual flow predicted by DNN surrogate (abscissa) 
are close to the surrogate boundary (ordinate).
Moreover, many of candidates are outside true feasible region examined by the simulator, i.e., these candidates are falsely classified as feasible ones by neural network (see yellow dots in Fig. \ref{fig:6para-optimize-boundary-divergence}).
Therefore these candidates are close to the true boundary.
For such a problem, obtaining an accurate surrogate classifier is key to our optimization. Therefore, our iterative training algorithm, which can adaptively improve the accuracy of the DNN classifier, is a key procedure for a good performance of our DiDo approach.

\section{Toy Example: Harmonic function}

To verify the validity of our method in solving high-dimensional data-informed optimization problem. Inspired by the practical problem of the rotor design, we construct a 100-dimensional optimization problem, whose optimal points locates on the boundary of the feasible region. 

\subsection{Problem description}
We consider an optimization problem, where the objective function $f(\vx)$ is a harmonic function $f(\vx)=-(x_1^2-\frac{1}{d-1}\sum_{i=2}^{d}x_{i}^2)$, which satisfies Poisson's equation $\nabla^2 f(\vx)=0$. The feasible region is 
$\Omega=\{\vx|||\vx||\leq 1\}=\{\vx=(x_1,...,x_d)^T|\sqrt{\sum_{i=1}^{d}{x_i^2}}\leq 1\}$.

The toy optimization problem is as follows, 
\begin{align}
&\min_{\vx}\quad -(x_1^2-\frac{1}{d-1}\sum_{i=2}^{d}x_{i}^2) \notag\\
&\hspace{0.1cm}\mathrm{s.t.}\hspace{0.3cm}\sqrt{\sum_{i=1}^{d}{x_i^2}}\leq 1
\label{subprob3}
\end{align}
where $\vx=(x_1,...,x_{d})^T\in \mathbb{R}^{d}$. For demonstration, we take $d=100$.

Note that the  harmonic function $f(\vx)$ satisfies extremum principle, which indicates that the minimum of problem ($\ref{subprob3}$) is achieved on the boundary.
As for the given case, it is clear that the minimum $-1$ is obtained at $\vx=(1,0,...,0)^T$ and $\vx=(-1,0,...,0)^T$. Remark that although the objective function and the constraints are analytically known, we assume that the objective function and the constraint functions can only be measured through sampling.

\subsection{Feasible region learned by a DNN classifier}


Similarly to the rotor problem, with the same settings, we first train a DNN classifier to learn the feasible region. 
We set initial sample size $n_0=3000$,initial sample region $B=[-0.173,0.173]^{100}$ and  we set $n_1=5000$ samples in each iteration. 
By algorithm \ref{alg:iterative method}, we obtain a well-trained classifier $f_{\vtheta_{\mathrm{c}}}(\vx) \in [0,1]$.
We  use the surrogate feasible region  $\{\vx|f_{{\vtheta_{\mathrm{c}}}}(\vx)\geq 0.5\}$ to represent the true feasible region $\Omega$.

\begin{figure}[!htb]
    \centering
    \subfigure[]{
    \includegraphics[scale=0.13]{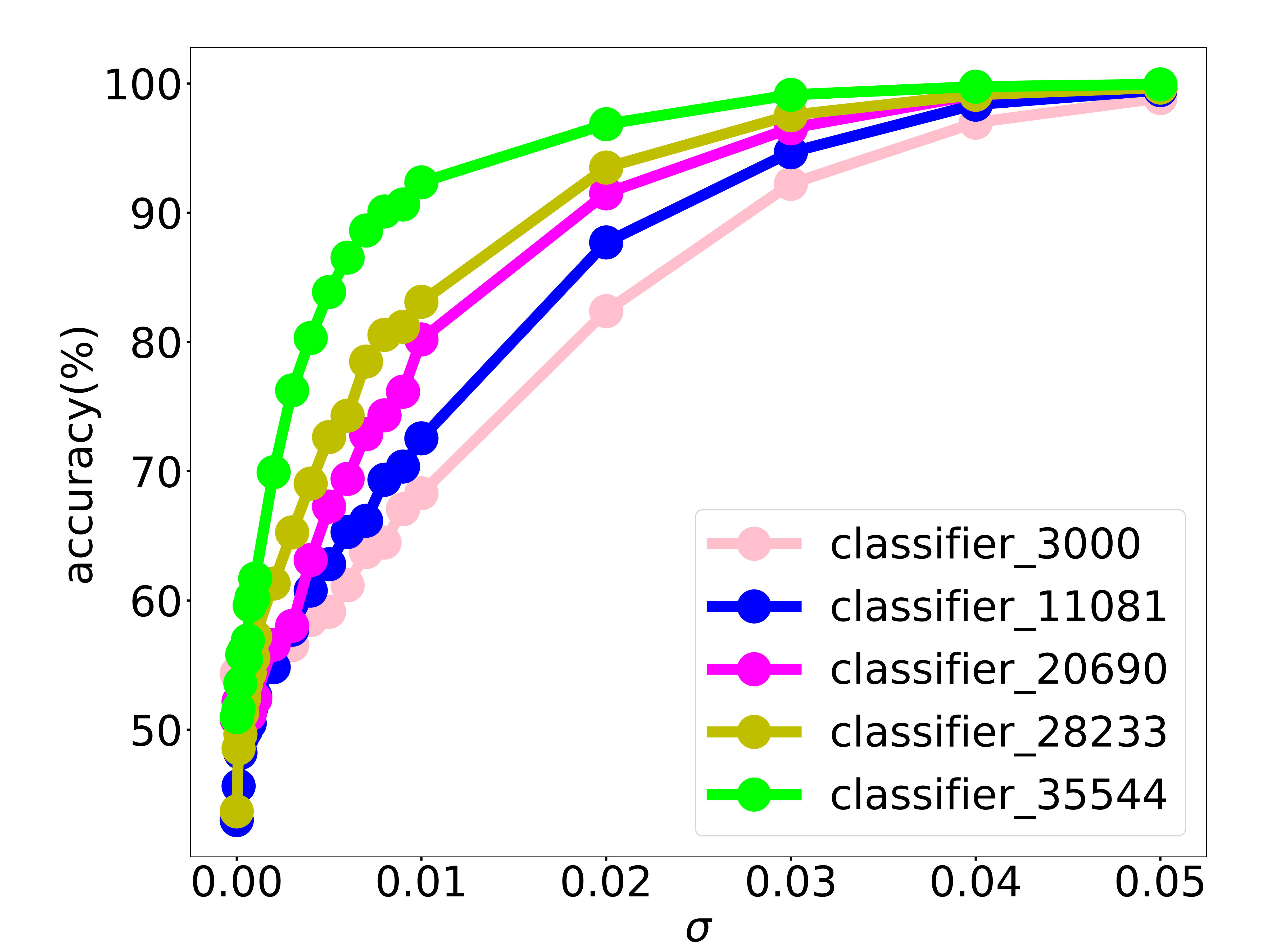}
    }
    \subfigure[]{
    \includegraphics[scale=0.13]{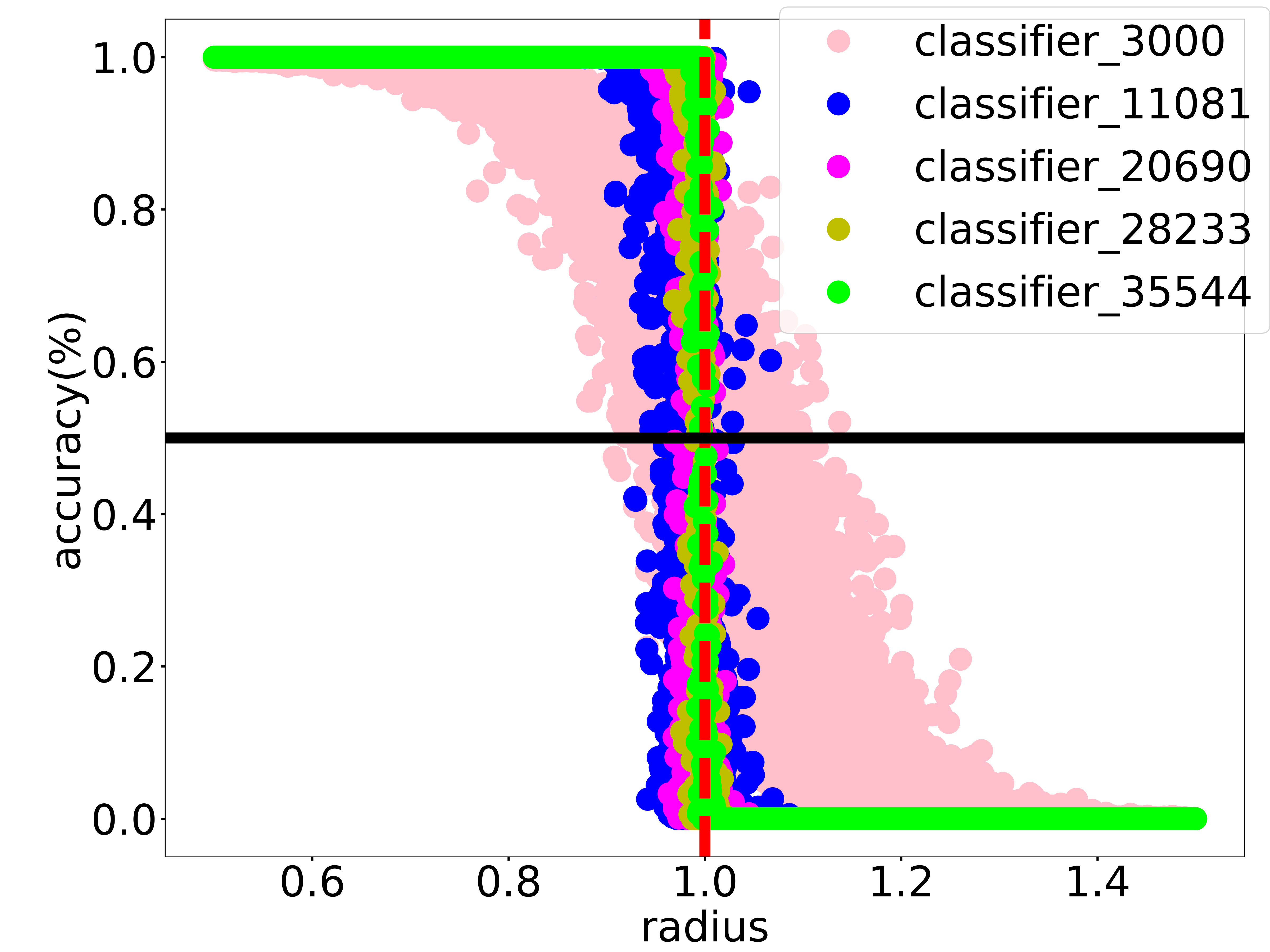}
    }

    \caption{The performance of the classifier. (a) Classification accuracy of the DNN classifier on the perturbed terms during iteration. Note that, there are not all iteration results and at each iteration $t$, we apply an extra constraint  $|f_{\vtheta_{\mathrm{c}}}^{(t)}(\vx_i)-0.5|\leq 0.1\}$ to the points sampled by LMC. In the two figures, label accuracy means classification accuracy after perturbation. As we add more data, the magnitude of the perturbed term when classifier accuracy on perturbed term increase from $50\%$ sharply gets smaller, which means the distance between the true boundary and surrogate boundary gets smaller,i.e.,  the performance of classifier is better; (b) The classifier values on the points uniformly distributed along the radial direction. As the iteration proceeds, the classifier is more closed to the real classification function $I(r\leq1)$.}
    \label{fig:100d-classifier-performance}
\end{figure}

During the iterative training, the accuracy of the classifier with Gaussian noise perturbation efficiently improves as shown in Fig. \ref{fig:100d-classifier-performance}(a).
In addition, for this toy example, we know the real feasible region is a unit ball and it is clear to visualize the boundary along the radial direction. Thus, we calculate $f_{\vtheta_{\mathrm{c}}}(r\vx)$, where $\vx$ is uniformly sampled on the real boundary of the feasible region and $r$ follows uniform distribution on the interval [0,2]. As shown in Fig. \ref{fig:100d-classifier-performance}(b), throughout the iterative training, the surrogate classifier approximates the true feasible region $I(r\leq1)$ better and better.

\subsection{Objective function learned by a DNN model}

\begin{figure}[!htb]
    \centering
    \includegraphics[scale=0.15]{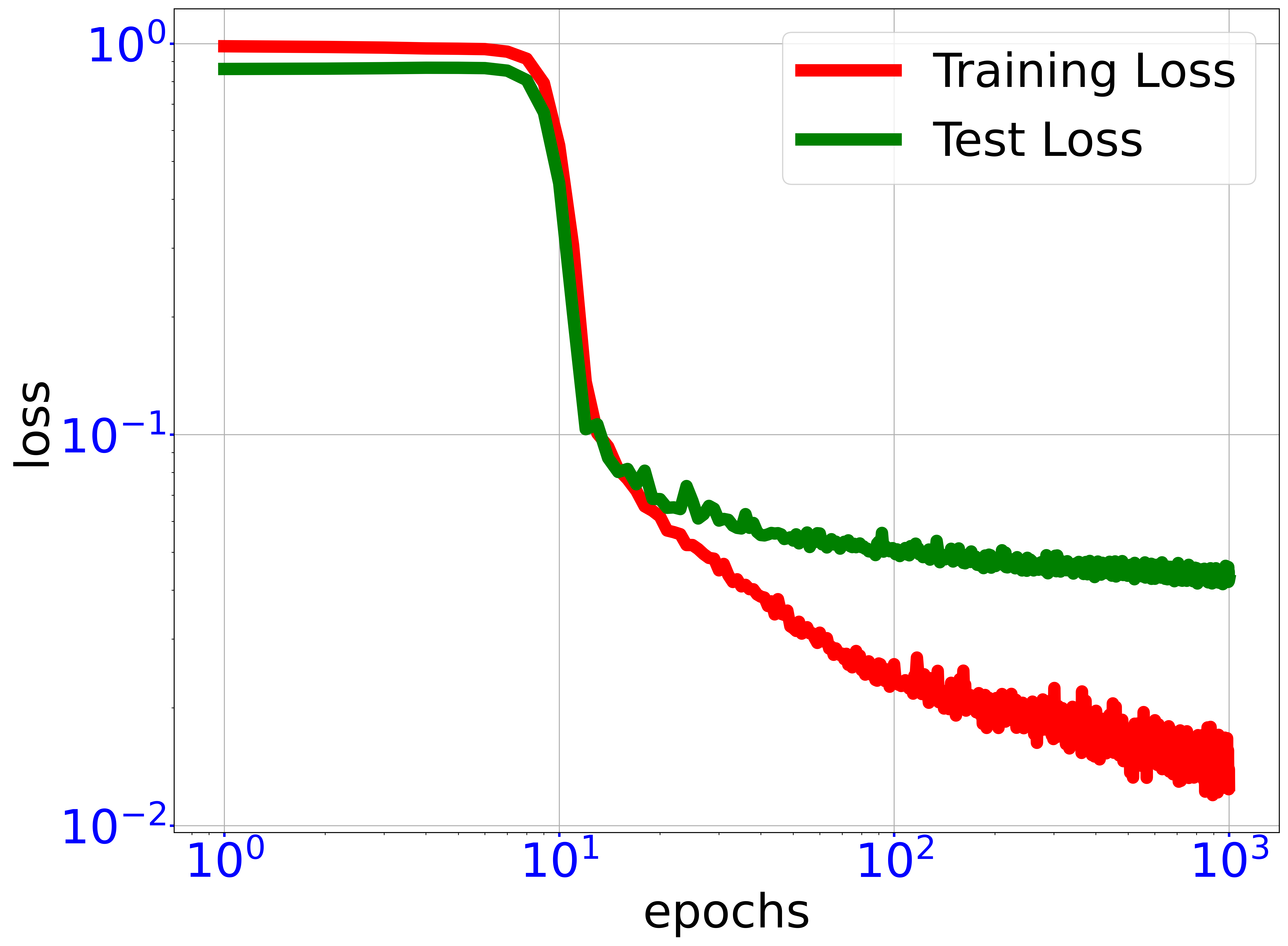}
    \caption{Training loss and test loss during training  DNN surrogate $f_{\vtheta_{\mathrm{o}}}(\vx)$ of toy example.}
    \label{fig:toyDNNloss}
\end{figure}

 By algorithm \ref{alg:sample method}, we use  the classifier $f_{\vtheta_{\mathrm{c}}}(\vx)$ obtained above to generate a training set $D_{\mathrm {obj}}$ of size $n_{t'} = 5,000$  
and train a GELU-DNN $f_{\vtheta_{\mathrm{o}}}(\vx)$ of hidden layer size 2000-1000-600-400-200. 
The test accuracy of the DNN $f_{\vtheta_{\mathrm{o}}}(\vx)$  is evaluated on a test set consisting of 2000 samples. 
As shown in  Fig. \ref{fig:toyDNNloss}, after training, the normalized RMSE training error is  $\sim 0.01$ whereas the normalized RMSE test error is $\sim 0.04$. 

\subsection{Deep optimization}
\begin{figure}[!htb]
    \centering
    \includegraphics[scale=0.15]{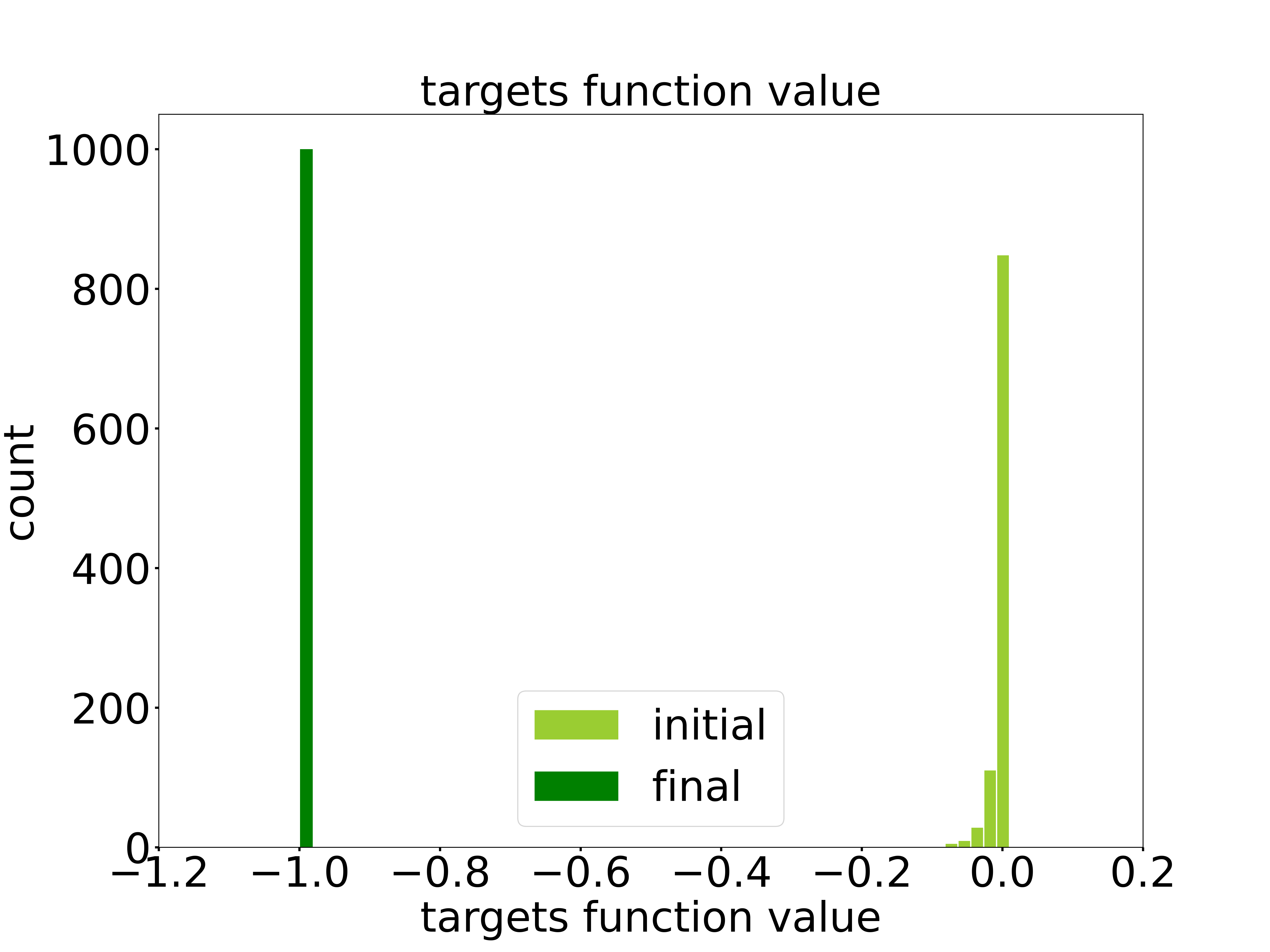}
    \caption{The distribution of the objective function values. Comparison between the objective function values on the initial points, i.e., the training samples used for learning DNN surrogate, and that on final candidates of optimal parameters.}
    \label{fig:laplace-optimization-result}
\end{figure}

With DNN surrogate $f_{\vtheta_{\mathrm{o}}}(\vx)$ and the well-trained DNN classifier $f_{\vtheta_{\mathrm{c}}}(\vx)$, by algorithm \ref{alg:Deepoptimization}, we obtain a set of candidates from different initial points. Note that we set the training samples used for learning DNN surrogate as initial points. For visualization, in Fig. \ref{fig:laplace-optimization-result}(a), we show the distribution of the objective function values of the initial points as well as that of the corresponding  candidates of optimal parameters. Note that the minimum objective function value among training samples used for learning DNN surrogate $\sim-0.1$, whereas
the objective function values of the  candidates of optimal parameters concentrate around $-0.98$ very close to the true minimum $-1$ of this problem.

\section{Conclusion and discussion}

Data-informed optimization problems are common in science or industry. Though intuitive, our DiDo approach provides a promising way to solve high-dimensional data-informed optimization problems. 
It is hoped that this work can help researchers and engineers in solving diverse types of high dimensional optimization problems where objective function and feasible region are implicit. 

For  a type of high dimensional optimization problems, whose optimal points located in the interior region, e.g., maximize the Gaussian function in a unit ball, we find that it is more difficult to sample sufficient useful points to fit the objective function well.
This phenomenon is due to the  concentration phenomena in high dimension space \cite{blum2020foundations}.  For example, if we uniformly sample data in an unit ball, the samples concentrate at an O(1/d) shell of the surface.
In practice, this phenomenon can be alleviated by using a proper sampling distribution, say radial uniform sampling, according to prior knowledge.

In the iterative training process, the hyperparameter $\beta$ in LMC is important to sample diverse points close to the surrogate boundary. If $\beta$ is too large, we observe that the added points  
concentrate at the surrogate decision boundary and the new classifier can even become less accurate.
This phenomenon is related to frequency principle, i.e.,  the  points close to the boundary are high frequency
in nature
, thus may result in worse generalization performance \cite{xu2020frequency}. 
Empirically, proper $\beta$ is needed for a steady improvement of accuracy of the DNN classifier.


\section*{Acknowledgments}
This work is sponsored by the National Key R\&D Program of China  Grant No. 2019YFA0709503 (Z. X.), the Shanghai Sailing Program, the Natural Science Foundation of Shanghai Grant No. 20ZR1429000  (Z. X.), the National Natural Science Foundation of China Grant No. 62002221 (Z. X.), Shanghai Municipal of Science and Technology Project Grant No. 20JC1419500 (Y. Z.), 
Shanghai Municipal of Science and Technology Major Project No. 2021SHZDZX0102, and the HPC of School of Mathematical Sciences and the Student Innovation Center at Shanghai Jiao Tong University.

\bibliographystyle{plain}
\bibliography{myRef}

\end{document}